\documentclass[a4paper]{amsart}
\usepackage{graphicx}
\usepackage{amsmath,amssymb,amscd}

\newtheorem{theorem}{Theorem}
\newtheorem{corollary}[theorem]{Corollary}

\newtheorem{lemma}[theorem]{Lemma}

\def\rank{\mathop{\rm rank}\nolimits}

\def\Aut{\mathop{\rm Aut}\nolimits}
\def\id{\mathop{\rm Id}\nolimits}
\def\diag{\mathop{\rm diag}\nolimits}
\def\gcd{\mathop{\rm gcd}\nolimits}

\title{A Canonical Form for a Symplectic Involution}

\author{H.W. Braden}

\address{
Maxwell Institute and School of Mathematics, The University of Edinburgh, 
James Clerk Maxwell Building, Peter Guthrie Tait Road, 
Edinburgh EH9 3FD, United Kingdom
}
\email{hwb@ed.ac.uk}

\begin{document}
\maketitle
\begin{abstract}
We present a canonical form for a symplectic involution $S\in Sp(2g,\mathbb{Z})$, $S^2=\id$;
the construction is algorithmic. Application is made in the Riemann surface setting.
\end{abstract}

\section{Introduction}

Canonical forms for matrices with integer coefficients are useful
in many settings: one may think of the Smith Normal Form or
Frobenius's decomposition of a skew matrix \cite{new72}, both of
which will be used later in the paper. In his study of real
abelian varieties Comessatti \cite{com25, com26} introduced a
canonical form for an involution in $GL(n,\mathbb{Z})$ (the
precise result will be recalled later). Here we shall establish
the symplectic analogue of Comessatti's theorem providing a
canonical form for a symplectic involution, $S\in
Sp(2g,\mathbb{Z})=\{\gamma\in Gl(2g,\mathbb{Z})\,\big|\,
\gamma\sp{T} J\gamma=J\}$, where throughout $J=
\begin{pmatrix}0&1_g\\-1_g&0\end{pmatrix}$ is the canonical symplectic pairing.
The canonical form with an immediate corollary is given by:
\begin{theorem}\label{sympcf} Let $S\in Sp(2g,\mathbb{Z})$ be a symplectic involution,
$S\sp{T}JS=J$ and
$S\sp2=\id$. Then $S$ is symplectically
equivalent to one of the form $S=\begin{pmatrix}a&0\\
0&a\end{pmatrix}$ where
\begin{equation}\label{canform}
a=\begin{pmatrix}1_{p}\\
& -1_{m}\\
&&\begin{matrix}0&1\\
1&0\end{matrix}\\
&&&\ddots\\
&&&&\begin{matrix}0&1\\
1&0\end{matrix}\\
\end{pmatrix}.
\end{equation}
If $t$ is the number of $2\times2$ blocks then $g=p+m+2t$.
\end{theorem}
\begin{corollary}\label{symplag} Let $S$ be a symplectic involution of
$W=\mathbb{Z}\sp{2g}$ with canonical pairing. Then $W=L_1\oplus
L_2$, $<L_i,L_i>=0$, with stable Lagrangian subspaces $SL_i=L_i$.
\end{corollary}
Some special cases of Theorem \ref{sympcf} are known in the
context of Riemann surfaces. The canonical form yields a different
proof of
\begin{theorem} Let $S$ be a conformal involution of the Riemann surface
${\mathcal{C}}$ of genus $g$ with $k$ fixed points and let
$\mathcal{C}'={\mathcal{C}}/<S>$ be the quotient surface of genus
$g'$. We may find a homology basis for ${\mathcal{C}}$ in which
$S$ takes the form (\ref{canform}) where $ g=p+m+2t$, $g'=p+t$ and
either
\begin{enumerate}
\item $k>0$ and $p=0$ whence $g'=t$ and $m=k/2-1$,
\item  $k=0$, $p=1$ whence $g'=t+1$ and $m=0$.
\end{enumerate}
\end{theorem}
The former case yields a result of Gilman \cite{gil73} while the
latter yields Fay's example of an unbranched covers \cite[Ch
IV]{fay73}. 

We note that the proof we present of Theorem \ref{sympcf} is constructive.
Before turning to the proofs we will give some further background including Comessatti's result that we will employ. If one could find a module-theoretic proof of Corollary \ref{symplag} then Theorem \ref{sympcf}
would follow from Comessatti's theorem.

Finally we remark that the theorems described in this paper are of significant utility in the computational
study of Riemann surfaces. Relevant for this volume Edge often studied curves and geometric configurations with high symmetry such as Klein's curve \cite{e44}, Bring's curve \cite{e78}
and the Fricke-Macbeath curve \cite{e67}; adapted homology bases may be found using the theorems outlined in this paper for curves with symmetries that, for example, significantly simplify the period matrices
and the calculation of the vector of Riemann constants \cite{bn10, bn12}.

\section{Background}
In order to place the result in context its helpful to see the
parallel between several results for the general linear and
symplectic groups. First,
\begin{lemma}\label{gcdunz}
If $\gcd(x_1,\ldots,x_m)=d$ then there exists $U\in
GL(m,\mathbb{Z})$ such that
$$(x_1,\ldots,x_m)U=(d,0,\ldots,0).$$
\end{lemma}
\begin{lemma}\label{gcdspnz}
If $\gcd(x_{1},x_{2},\ldots,x_{2g})=d$ then there exists $S\in
Sp(2g,\mathbb{Z})$ such that
$$(x_{1},x_{2},\ldots,x_{2g})S=(d,0,\ldots,0).$$
\end{lemma}
Lemma \ref{gcdunz} is classical (see for example \cite{new72}).
Lemma \ref{gcdspnz} seems less well-known; the first proof of this
I am aware of is \cite{rei54}.

In his study of real abelian varieties Comessatti \cite{com25,
com26} introduced the following canonical form.
\begin{theorem}\label{commthm}Let $M$ be a free $\mathbb{Z}$-module of rank $m$ and let
$S\in\Aut(M)$ be an involution. Let $\rank
M_\pm:=s_\pm$ where $M_\pm$ are the submodules
\begin{equation*}
\label{submodS}
M_+:=\{x\in M\,|\ Sx=x\},\qquad
M_-:=\{x\in M\,|\ Sx=-x\}.
\end{equation*}
Then we may find a basis of $M$
such that $S$ takes the form
\begin{equation}\label{comessati2}
S=\begin{pmatrix}1_{s_+ -\lambda}&\\
& -1_{s_- -\lambda}\\
&&\begin{matrix}0&1\\
1&0\end{matrix}\\
&&&\ddots\\
&&&&\begin{matrix}0&1\\
1&0\end{matrix}
\end{pmatrix}
\end{equation}
where we have $\lambda$ copies of the matrix $\begin{pmatrix}0&1\\
1&0\end{pmatrix}$.  Moreover, $s_+$,
$s_-$ and $\lambda$ are invariants of $S$.
\end{theorem}
Here $\lambda$ is known as the Comessatti character. (A modern review of
Comessatti's work may be found in \cite{cp96}.) Silhol \cite{sil82} expressed these
invariants in terms of the group cohomology of
$G=<1,S>$. Then
\begin{equation*}\label{gpcohcomess}
H\sp{i}(G,M)=\begin{cases}
M_+ \cong(\mathbb{Z}_2)\sp{s\sp+} &i=0,\\
\dfrac{M_-}{(1-S)M}\cong(\mathbb{Z}_2)\sp{s\sp--\lambda}&i\equiv1 \pmod2,\\
\dfrac{M_+}{(1+S)M}\cong(\mathbb{Z}_2)\sp{s\sp+-\lambda}&i\equiv0 \pmod2,\ i>0.
\end{cases}
\end{equation*}

In the study of real structures the focus of attention are
anti-holomorphic involutions $S\sp{T}JS=-J$, $S^2=\id$ rather than
holomorphic involutions. Comessatti showed that an anti-holomorphic  involution $S$ takes
the form $\begin{pmatrix}1_{g}&H\\0&-1_{g}\end{pmatrix}$ where $H$
is a symmetric bilinear form over $\mathbb{Z}_2$. These forms are
determined by the rank of $H$ and whether $\diag(H)$ is nonzero or
not. We have either (see for example \cite{vin93})
$$H=\begin{pmatrix}0&1& \\ 1&0&\\ &&\ddots&\\
&&&0&1&\\ &&&1&0&\\ &&&&&0&\\ &&&&&&\ddots&\\
&&&&&&&&0\end{pmatrix},\qquad
H=\begin{pmatrix}1& \\  &\ddots&\\
&&1&\\ &&&0&\\  &&&&&\ddots&\\
&&&&&&&0\end{pmatrix}.
$$
An algorithm that constructs such a basis for a Riemann surface with real structure is given in
\cite{kk14}.

\emph{Remark:} Comessatti's theorem admits both purely
module-theoretic and constructive proofs. The proof that follows
of Theorem \ref{sympcf} is constructive. If one could find a
module-theoretic proof of Corollary \ref{symplag} then the theorem
would follow from Comessatti's theorem.

\section{Proof of theorem 1.}
The proof is constructive. Writing $S=\begin{pmatrix}
a & b \\
c& d
\end{pmatrix}$
where $a$, $b$, $c$, $d$ are block $g\times g$ integer matrices,
the constraints $S\sp{T}JS=J$ and $S\sp2=\id$ mean that $S$ takes
the form
\begin{equation}\label{holinvnform}
S=\begin{pmatrix}
a & b \\
c& a\sp{T}
\end{pmatrix},\qquad 0=b+b\sp{T}=c+c\sp{T}=ab+ba\sp{T}=ca+a\sp{T}c,\ a^2+bc=\id.
\end{equation}
We remark that if $U\in GL(g,\mathbb{Z})$ and $\mu=\mu\sp{T}$ then
the rotations and translations
$$R_U=\begin{pmatrix} U&0\\ 0&U\sp{-1\, T}\end{pmatrix},
\qquad T_\mu=\begin{pmatrix} 1&\mu\\ 0&1
\end{pmatrix}$$ are symplectic. In particular
the similarity transformation
$$R_{U\sp{-1}}SR_U=\begin{pmatrix}
U\sp{-1}aU & U\sp{-1}bU\sp{-1\,T} \\
U\sp{T}cU& U\sp{T}a\sp{T}U\sp{-1\,T}
\end{pmatrix}
$$
yields a similarity transformation on $a$ and takes $c$ to a
congruent matrix.

The first step of the proof is to make a symplectic transformation
so that $c=0$. The Frobenius decomposition of the skew matrix $c$
says there exists $U\in GL(g,\mathbb{Z})$ such that $c=U\sp{T}DU$
where the congruent matrix $D$ takes the form
$$D=\begin{pmatrix}
  {0} & {d_{1}} & & & & & & & & \\
  {-d_{1}} & {0} & & & & & & & & \\
   & & {0} & d_{2} & & & & & & \\
   & & {-d_{2}} & {0} & & & & & & \\
   & & & & {\ddots} & {\ddots} & & & & \\
   & & & & {\ddots} & {\ddots}& & & & \\
   & & & & & & {0} & {d_{s}} & & \\
   & & & & & & {-d_{s}} & {0} & & \\
   & & & & & & & & & {0} \\
\end{pmatrix}
$$
with $d_{i}| d_{i + 1}$ for  $1 \leq i \leq s-1$  and $\rank\,
c=2s$. By using the appropriate symplectic transformation $R_U$ we
may suppose that $c$ is in the Frobenius form $D$ stated. Let
$\gcd(a_{21},d_1)=\nu$. Then there are $p,q,u,v\in\mathbb{Z}$ such
that
$$a_{21}=\nu p,\qquad d_1=\nu q,\qquad up-vq=1.$$
Then the symplectic matrix
$$T=\begin{pmatrix}
1&0&0&0&0&0\\ 0&u&0&0&v&0\\ 0&0&1_{g-2}&0&0&0\\
0&0&0&1&0&0\\ 0&q&0&0&p&0\\ 0&0&0&0&0&1_{g-2}
\end{pmatrix},$$
is such that
$$
T.S.T^{-1}=\begin{pmatrix}
a' & b' \\
D'& a'\sp{T}
\end{pmatrix},\qquad
D'=\begin{pmatrix}
  {0} & 0 & & & & & & & & \\
  0 & {0} & & & & & & & & \\
   & & {0} & d_{2} & & & & & & \\
   & & {-d_{2}} & {0} & & & & & & \\
   & & & & {\ddots} & {\ddots} & & & & \\
   & & & & {\ddots} & {\ddots}& & & & \\
   & & & & & & {0} & {d_{s}} & & \\
   & & & & & & {-d_{s}} & {0} & & \\
   & & & & & & & & & {0} \\
\end{pmatrix}
.$$ Continuing in this way we see that $S$ is similar via a
symplectic transformation to the case when $c=0$.

With $c=0$ we see from from (\ref{holinvnform}) that $a^2=1$.
Using the freedom to make a similarity transform to $a$, noted
above, we may now use Comessatti's theorem \ref{commthm} to put
$a$ into the canonical form
$$a
=\begin{pmatrix}1_{r}\\
& -1_{s}\\
&&\begin{matrix}0&1\\
1&0\end{matrix}\\
&&&\ddots\\
&&&&\begin{matrix}0&1\\
1&0\end{matrix}\\
\end{pmatrix}
$$
for appropriate $r$ and $s$. At this stage we have that
$$
S=\begin{pmatrix}
a & b \\
0& a\sp{T}
\end{pmatrix},\qquad 0=b+b\sp{T}=ab+ba\sp{T}
$$ and in block form
$$a=\begin{pmatrix}1_r&&\\
& -1_s\\ &&Q\end{pmatrix}, \qquad
Q=\begin{pmatrix}\begin{matrix}0&1\\
1&0\end{matrix}\\
&\ddots\\
&&\begin{matrix}0&1\\
1&0\end{matrix}\\
\end{pmatrix},
$$
where $Q$ is a $2l\times 2l$ matrix. Now solving for
$0=b+b\sp{T}=ab+ba\sp{T}$ we find that $b$ has the form
$$
b=\begin{pmatrix}0&x&y\\ -x\sp{T}&0&z\\ -y\sp{T}&-z\sp{T}
&\gamma\end{pmatrix},\qquad
\gamma+\gamma\sp{T}=0=\gamma Q+Q\gamma,\ \ y=-yQ,\ z=zQ.
$$
Here $x\in M_{r,s}(\mathbb{Z})$, $y\in M_{r,2l}(\mathbb{Z})$,
$z\in M_{s,2l}(\mathbb{Z})$, $\gamma\in M_{2l,2l}(\mathbb{Z})$.
Thus each row of the matrix $y$ takes the form
$$(y_{i1},-y_{i1},y_{i2},-y_{i2},\ldots,y_{il},-y_{il}),\qquad
1\le i\le r,$$ while each row of the matrix $z$ takes
the form
$$(z_{j1},z_{j1},z_{j2},z_{j2},\ldots,z_{jl},z_{jl}),\qquad
1\le j\le s.$$ Further the skew-symmetric matrix
$\gamma$ may be written as $2\times2$ blocks
$$\gamma=\begin{pmatrix}d_{11}&m_{12}&\ldots&m_{1l}\\
-m_{12}\sp{T}&d_{22}&\\
\vdots&&\ddots\\
-m_{1l}\sp{T}&&&d_{ll}\\
\end{pmatrix}, \
d_{ii}=\begin{pmatrix}0&\alpha_i\\
-\alpha_i&0\end{pmatrix},\ m_{ij}=\begin{pmatrix}\beta_{ij}&\delta_{ij}\\
-\delta_{ij}&-\beta_{ij}\end{pmatrix}.
$$

Observe that
$$T_\mu S T_{-\mu}=\begin{pmatrix}
a & b+\mu a\sp{T}-a \mu \\
0& a\sp{T}
\end{pmatrix}
$$
and so if $$\mu=\begin{pmatrix} \mu_1&\mu_2&\mu_3\\
\mu_2\sp{T}&\mu_4&\mu_5\\ \mu_3\sp{T}&\mu_5\sp{T}&\mu_6
\end{pmatrix},\qquad \mu_1=\mu_1\sp{T},\ \mu_4=\mu_4\sp{T}, \mu_6=\mu_6\sp{T}
$$
then $\mu=\mu\sp{T}$ and
$$\mu a\sp{T}-a \mu =
\begin{pmatrix} 0&-2\mu_2&\mu_3Q-\mu_3\\
2\mu_2\sp{T}&0&\mu_5 Q+\mu_5\\ \mu_3\sp{T}-Q
\mu_3\sp{T}&-\mu_5\sp{T}-Q\mu_5\sp{T}&\mu_6 Q-Q\mu_6
\end{pmatrix}.
$$
Thus if we choose the rows of the matrix $\mu_3$ to be
$(y_{i1},0,y_{i2},0,\ldots,y_{il},0)$ ($
1\le i\le r$) then
$$y+\mu_3Q-\mu_3=0.$$ Similarly if
the rows of the matrix of the matrix $\mu_5$ to be
$-(z_{j1},0,z_{j2},0,\ldots,z_{jl},0)$ ($ 1\le j\le s$) then
$$z+\mu_5 Q+\mu_5=0.$$
Finally taking $\mu_6$ to be of the form
$$\mu_6=\begin{pmatrix}d'_{11}&m'_{12}&\ldots&m'_{1l}\\
{m'}_{12}\sp{T}&d'_{22}&\\
\vdots&&\ddots\\
{m'}_{1l}\sp{T}&&&d'_{ll}\\
\end{pmatrix}, \
d'_{ii}=\begin{pmatrix}-\alpha_i&0\\
0&0\end{pmatrix},\ m'_{ij}=\begin{pmatrix}-\delta_{ij}&-\beta_{ij}\\
0&0\end{pmatrix}.
$$
yields
$$\gamma+\mu_6 Q-Q\mu_6=0.$$
Therefore we may take $b$ to be of the form form
$$
b=\begin{pmatrix}0&x&0\\ -x\sp{T}&0&0\\ 0&0
&0\end{pmatrix},
$$
and where $x$ is a $(0,1$)-matrix.

At this stage we have shown that we may choose a symplectic basis
in which the involution $S$ takes the block form
$$S=\begin{pmatrix}1_r&0&0&0&x&0 \\
0& -1_s&0&-x\sp{T}&0&0\\ 0&0&Q&0&0&0\\
0&0&0&1_r&0&0\\
0&0&0&0& -1_s&0\\
0&0&0&0&0&Q\end{pmatrix}
$$
where $x$ is a $(0,1$)-matrix. Further, by use of the rotation
$R_U$ with $U$ of the form
$$U\sp{-1}=\begin{pmatrix}A&&\\
& B\\ &&1\end{pmatrix}, \qquad A\in GL(r,\mathbb{Z}),\ B\in
GL(s,\mathbb{Z}),$$ we may transform $x$ to $AxB\sp{T}$. Making use
of the Smith normal form and the ability to remove even integral
parts of $x$ by a translation we may therefore assume $x$ to have
only $1$'s and $0$'s along the diagonal and be zero off the
diagonal. Suppose there are $t\le \min(r,s)$ $1$'s on the
diagonal. Then we may write
$$
S=\begin{pmatrix}1_{r-t}&0&0&&& \\
0& -1_{s-t}&0&&&\\ 0&0&Q&\\
&&&1_{r-t}&0&0\\
&&&0& -1_{s-t}&0\\
&&&0&0&Q\end{pmatrix}\oplus S'\oplus \ldots\oplus S'
$$
where we have $t$ copies of the symplectic matrix $$
S'=\begin{pmatrix}1&0&0&1 \\
0& -1&-1&0\\
0&0&1&0\\
0&0&0& -1\end{pmatrix}
$$
and we are indicating a symplectic decomposition. Now consider
$$V=\begin{pmatrix}1&0&0&0\\
1&0&0&1 \\
0& 1&1&0\\
0&-1&0&0\end{pmatrix}.$$ Then
$$V\sp{T}
\begin{pmatrix}0&1_2\\-1_2&0\end{pmatrix}V=
\begin{pmatrix}0&1_2\\-1_2&0\end{pmatrix},\qquad
V S' V^{-1}=\begin{pmatrix}0&1&0&0\\
1&0&0&0 \\
0& 0&0&1\\
0&0&1&0\end{pmatrix}.
$$
Thus by conjugation we may bring $S$ to the desired form and have
established the theorem.

\section{Proof of theorem 3.}
We now apply our results in the setting where we have a Riemann
surface ${\mathcal{C}}$ of genus $ g>0$ with nontrivial finite
group of symmetries $G\le\Aut {\mathcal{C}}$. ($\Aut
{\mathcal{C}}$ is necessarily finite for $g\ge2$.) $\Aut
{\mathcal{C}}$ acts naturally on ${\mathcal{C}}$, $H_1(
{\mathcal{C}},\mathbb{Z})$ and the harmonic differentials.
Consider the quotient Riemann surface
$\pi:{\mathcal{C}}\rightarrow\mathcal{C}'={\mathcal{C}}/G$ of
genus $g'$. From (\ref{canform}) $g=(p+t)+(m+t)$ and we can form
$p+t$ invariant differentials and $m+t$ anti-invariant
differentials under the action of $S$; then $g'=p+t$ is the genus
of $\mathcal{C}'$. By Riemann-Hurwitz if there are $k\ge0$ fixed
points of $S$ then $ g-1 =2(g'-1)+k/2$ yields
$$m=p+k/2-1.$$

Hurwitz showed that $\phi\in\Aut \mathcal{C}$ is the identity if
and only if it induces the identity on
$H_1({\mathcal{C}},\mathbb{Z})$. Accola \cite{acc67} strengthened
this result and showed that for $g\ge2$ if there exist two pairs
of canonical cycles such that (in homology)
$\phi({\mathfrak{a}}_1)={\mathfrak{a}}_1$,
$\phi({\mathfrak{a}}_2)={\mathfrak{a}}_2$,
$\phi({\mathfrak{b}}_1)={\mathfrak{b}}_1$ and
$\phi({\mathfrak{b}}_2)={\mathfrak{b}}_2$ then $\phi$ is the
identity. (Simpler proofs of this result were obtained by Earle as
well as Grothendieck and Serre, see \cite{gil77}.) Accola's result
means in the canonical form for the symplectic involution above we
have $p\le1$. We will have therefore proven the theorem once we
establish
\begin{lemma} If $k>0$ then $p=0$.
\end{lemma}
\emph{Proof of Lemma.} Let $\{\gamma_a\}_{a=1}\sp{2 g}$  be a
basis for $H_1({\mathcal{C}},\mathbb{Z})$ ordered such that
$\gamma_a=\mathfrak{a}_a$, $\gamma_{ g+a}=\mathfrak{b}_a
$
($a=1,\ldots,g$) are canonically paired,
$<\mathfrak{a}_a,\mathfrak{b}_b>=\delta_{ab}$, and the symplectic
form is ${J}_{ab}=<\gamma_a,\gamma_b>$. Let ${\alpha}_b$ denote a
basis of the harmonic forms paired with the homology cycles
$\gamma_a$ by $\int_{\gamma_a}{\alpha}_b=\delta_{ab}$. With the
metric on (complexified as necessary) one-forms $ (\alpha, \beta)
= \int_{{\mathcal{C}}} \alpha \wedge \ast \bar{\beta}$ then we
also have that
\begin{equation}
{J}_{ab}=(\ast{\alpha}_a,{\alpha}_b)=-({\alpha}_a,
\ast{\alpha}_b)=\int_{{\mathcal{C}}}{\alpha}_a\wedge {\alpha}_b,
\end{equation}
where $\ast$ is the Hodge star operator. If $u,v\in
H\sp{1}(\mathcal{C}',\mathbb{R})$ then $|G|\,(u,\ast
v)=(\pi\sp\ast u,\ast\pi\sp\ast v)$. Letting
$\{\gamma_i'\}_{i=1}\sp{2 g'}$, $\{\alpha_i'\}_{i=1}\sp{2 g'}$
denote the analogous quantities for $\mathcal{C}'$ we may write
$u=\sum_i u_i\,\alpha_i'$ and similarly for $v$.

Suppose that $p>0$. Then (upon possible relabelling) we have
$S({\mathfrak{a}}_1)={\mathfrak{a}}_1$,
$S({\mathfrak{b}}_1)={\mathfrak{b}}_1$ and
$S\sp\ast{\alpha}_1={\alpha}_1$, $S\sp\ast{\alpha}_{
g+1}={\alpha}_{g+1}$. Now $\pi\sp\ast u$ for $u\in
H\sp{1}(\mathcal{C}',\mathbb{R})$ span the invariant differentials
of $H\sp{1}({\mathcal{C}},\mathbb{R})$ and we may find $u$, $v$
such that $\pi\sp\ast u={\alpha}_1$, $\pi\sp\ast v={\alpha}_{
g+1}$. We have that
\begin{equation}\label{liftip}
2\,(u,\ast v)={\rm{ord}}(S)\,(u,\ast v)= (\pi\sp\ast
u,\ast\pi\sp\ast v)=({\alpha}_{1},\ast{\alpha}_{
g+1})=-1.
\end{equation}

Now suppose that in addition there exists a fixed point $P$ of
$S$. Thus for all $Q$,
$$\int_P\sp{Q}{\alpha}_1=\int_P\sp{Q}\pi\sp\ast{\alpha}_1=
\int_P\sp{S(Q)}{\alpha}_1$$ and so consequently (as the path
between $Q$ and $S(Q)$ may be arbitrary)
$\int_Q\sp{S(Q)}{\alpha}_1\in \mathbb{Z}$. But now if $\gamma'$ is
any cycle on $\mathcal{C}'$ containing the arbitrary point
$\pi(Q)$ this may be lifted to a path in ${\mathcal{C}}$ beginning
at $Q$ and ending at $S\sp{l}(Q)$ for some $l$. (We may assume
this lifted path does not pass through any of the fixed points of
$S$.) Then
$$\int_{\gamma'} u=\int_Q\sp{S\sp{l}(Q)} \pi\sp\ast u=
\int_Q\sp{S\sp{l}(Q)}{\alpha}_1\in\mathbb{Z}.
$$
As this is true for any path $\gamma'$ we have $u=\sum_i
n_i\alpha_i'$ with $n_i\in \mathbb{Z}$ and similarly for $v=\sum_i
m_i\alpha_i'$ with $m_i\in \mathbb{Z}$. Therefore $(u,\ast
v)=-n\sp{T}Jm\in\mathbb{Z}$, but from (\ref{liftip}) this is not
possible. Thus if $p>0$ then $k=0$.

\hfill$\Box$

\providecommand{\bysame}{\leavevmode\hbox
to3em{\hrulefill}\thinspace}
\bibliographystyle{amsalpha}

\end{document}